\documentclass{article}
\usepackage{amssymb}
\usepackage{amsbsy}
\begin{document}
\centerline{\Large \bf } \vskip 6pt

\begin{center}{\Large \bf Orlicz hormonic Blaschke addition\footnote{Research is supported by
National Natural Science Foundation of China
(11371334).}}\end{center}

\vskip 10pt
\begin{center}
\centerline{Chang-Jian Zhao} \centerline{\it Department of
Mathematics, China Jiliang University, Hangzhou 310018, P. R.
China}\centerline{\it Email: chjzhao@163.com}
\end{center}

\vskip 10pt

\begin{center}
\begin{minipage}{12cm}
{\bf Abstract}~ Recently, Gardner, Hug and Weil have introduced
the Orlicz-Brunn-Minkowski theory: a general framework, additions,
and inequalities. Following this, in the paper we consider Orlicz
dual Brunn-Minkowski theory. We introduce {\it Orlicz hormonic
Blaschke addition} which is an extension of the $L_{p}$ hormonic
Blaschke addition and $L_{p}$ radial Minkowski addition,
respectively. Inequalities of dual Minkowski and Brunn-Minkowski
type are obtained for the {\it Orlicz hormonic Blaschke addition}.
The new Orlicz dual Brunn-Minkowski inequality implies the dual
and $L_{p}$-dual Brunn-Minkowski inequalities, respectively. New
Orlicz dual Minkowski inequality implies the $L_{p}$-dual
Minkowski inequality. One of these has connections with the
conjectured log-Brunn-Minkowski inequality of Lutwak, Yang, and
Zhang, and in fact we show a log dual Minkowski inequality.
Finally, we introduce the concept of Orlicz dual projection body
and an inequality similar to Orlicz projection body is
established.


\end{minipage}
\end{center}
\vskip 20pt

\vskip 10pt \noindent{\large \bf 1 ~Introducation}\vskip 10pt

The setting for this paper is $n$-dimensional Euclidean space
${\Bbb R}^{n}$. Let ${\cal K}^{n}$ denote the set of convex bodies
(compact, convex subsets with non-empty interior) in ${\Bbb
R}^{n}$ containing the origin in their interiors. We reserve the
letter $u$ for unit vectors, and the letter $B$ for the unit ball
centered at the origin. The surface of $B$ is $S^{n-1}$.

Throughout the paper, the standard orthonormal basis for ${\Bbb
R}^{n}$ will be $\{e_{1},\ldots,e_{n}\}$. Let $\Phi_{n},$
$n\in{\Bbb N}$, denote the set of convex functions
$\varphi:[0,\infty)^{n} \rightarrow [0,\infty)$ that are strictly
increasing in each variable and satisfy $\varphi(0)=0$ and
$\varphi(e_{j})=1$, $j=1,\ldots,n$. When $n=1$, we shall write
$\Phi$ instead of $\Phi_{1}$. The left derivative and right
derivative of a real-valued function $f$ are denoted by $(f)'_{l}$
and $(f)'_{r}$, respectively. We use $V(K)$ for the
$n$-dimensional volume of convex body $K$. Let $h(K,\cdot):
S^{n-1} \rightarrow {\Bbb R},$ denote the support function of a
convex body $K$; i.e. for $u\in S^{n-1}$, $h(K,u)=\max\{u\cdot x:
x\in K\},$ where $u\cdot x$ denotes the usual inner product $u$
and $x$ in ${\Bbb R}^{n}$. Let $\delta$ denote the Hausdorff
metric on ${\cal K}^{n}$, i.e., for $K, L\in {\cal K}^{n} ,$
$\delta(K,L)=|h_K-h_L|_{\infty},$ where $|\cdot|_{\infty}$ denotes
the sup-norm on the space of continuous functions $C(S^{n-1}).$

Associated with a compact subset $K$ of ${\Bbb R}^n$, which is
star-shaped with respect to the origin, is its radial function
$\rho(K,\cdot): S^{n-1}\rightarrow [0,\infty),$ defined for $u\in
S^{n-1}$, by $\rho(K,u)=\max\{\lambda\geq 0: \lambda u\in K\}.$ If
$\rho(K,\cdot)$ is positive and continuous, $K$ will be called a
star body. Let ${\cal S}_{o}^{n}$ denote the set of star bodies in
${\Bbb R}^{n}$. Let $\tilde{\delta}$ denote the radial Hausdorff
metric, as follows, if $K, L\in {\cal S}_{o}^{n}$, then
$\tilde{\delta}(K,L)=|\rho_{K}-\rho_{L}|_{\infty}$ (see e.g., [6]
and [40]).

In the $L_{p}$-Brunn-Minkowski theory, $L_{p}$ addition,
introduced by Firey in the 1960's. Denoted by $+_{p}$, this is
defined for $1\leq p\leq\infty$ by
$$h(K+_{p}L,x)^{p}=h(K,x)^{p}+h(L,x)^{p},\eqno(1.1)$$
for all $x\in {\Bbb R}^{n}$ and compact convex sets $K$ and $L$ in
${\Bbb R}^{n}$ containing the origin, where the functions are the
support functions of the sets involved (see e.g., [24]). When
$p=\infty$, (1.1) is interpreted as
$h(K+_{\infty}L,x)=\max\{h(K,x), h(L,x)\}$, as is customary. When
$p=1$, (1.1) defines ordinary Minkowski addition and then $K$ and
$L$ need not contain the origin.

$L_{p}$ addition and the $L_{p}$ Minkowski, Brunn-Minkwski
inequalities are fundamental inequalities from the $L_{p}$
Brunn-Minkowski theory. For recent important results and more
information from this theory, we refer to [14], [15], [16], [17],
[23], [25], [29], [30], [31], [32], [33], [36], [37], [43], [44],
[45] and the references therein.

In recent years, a new extension of $L_{p}$-Brunn-Minkowski theory
is Orlicz-Brunn-Minkowski theory, initiated by Lutwak, Yang, and
Zhang [34] and [35]. In these papers the fundamental notions of
$L_{p}$-centroid body and $L_{p}$-projection body were extended to
an Orlicz setting. It represents a generalization of the
$L_{p}$-Brunn-Minkowski theory, analogous to the way that Orlicz
spaces generalize $L_{p}$ spaces. The Orlicz centroid inequality
for star bodies was established in [48] which is an extension from
convex to star bodies. The other articles advance the theory can
be found in literatures [13], [19], [21] and [39].

Recently, Gardner, Hug and Weil ([8]) constructed a general
framework for the Orlicz-Brunn-Minkowski theory, and made clear
for the first time the relation to Orlicz spaces and norms. In
particularly, they defined the Orlicz sum $K+_{\varphi}L$ of
compact convex sets $K$ and $L$ in ${\Bbb R}^{n}$ containing the
origin, implicitly, by
$$\varphi\left(\frac{h(K,x)}{h(K+_{\varphi}L,x)},\frac{h(L,x)}
{h(K+_{\varphi}L,x)}\right)=1,\eqno(1.2)$$ for $x\in {\Bbb
R}^{n}$, if $h(K,x)+h(L,x)>0$, and by $h(K+_{\varphi}L,x)=0$, if
$h(K,x)=h(L,x)=0$. Here $\varphi\in \Phi_{2}$, the set of convex
functions $\varphi:[0,\infty)^{2}\rightarrow [0,\infty)$ that are
increasing in each variable and satisfy $\varphi(0,0)=0$ and
$\varphi(1,0)=\varphi(0, 1)=1$. Orlicz addition reduces to $L_{p}$
addition, $1\leq p<\infty$, when $\varphi(x_{1}, x_{2})=x_{1}^{p}
+x_{2}^{p}$, or $L_{\infty}$ addition, when $\varphi(x_{1},
x_{2})=\max\{x_{1}, x_{2}\}$. The particular instance of interest
corresponds to using (1.2) with $\varphi(x_{1},
x_{2})=\varphi_{1}(x_{1})+\varepsilon\varphi_{2}(x_{2})$ for
$\varepsilon>0$ and some $\varphi_{1},\varphi_{2}\in \Phi$, in
which case write $K+_{\varphi,\varepsilon}L$ instead of
$K+_{\varphi}L$.

Gardner, Hug and Weil [8] obtained the Orlicz-Brunn-Minkowksi
inequality for the Orlicz sum. In the end they also defined the
Orlicz mixed volume of convex bodies which contain the origin in
their interiors and got the Orlicz-Minkowski inequality for Orlicz
mixed volume and et al.

In the paper, we consider the same problem for dual theory of star
bodies not convex bodies by adopting a dual idea of Gardner, Hug
and Weil [8]. Let us first recall the concept, the harmonic
Blaschke addition, defined by Lutwak [26]. Suppose $K$ and $L$ are
star bodies in ${\Bbb R}^{n}$, the harmonic Blaschke linear
combination, $K\check{+}L,$ by
$$\frac{\rho(K\check{+}L,
\cdot)^{n+1}}{V(K\check{+}L)}=\frac{\rho(K,\cdot)^{n+1}}{V(K)}+\frac{\rho(L,\cdot)^{n+1}}{V(L)}.\eqno(1.3)$$
Lutwak [26] established the Brunn-Minkowski inequality for the
harmonic Blaschke addition. If $K,L\in {\cal S}_{o}^{n},$ then
$$V(K \check{+}L)^{1/n}\geq V(K)^{1/n}+V(L)^{1/n},
\eqno(1.4)$$ with equality if and only if $K$ and $L$ are dilates.
For a systematic investigation of fundamental characteristics of
additions of convex and star bodies, we refer to [7], [8] and [9].
The harmonic Blaschke addition and radial Minkowski addition are
the basis for the dual Brunn-Minkowski theory (see, e.g., [3],
[10], [11], [12], [20], [22], [26] and [41] for recent important
contributions).

Here, we first define the {\it Orlicz harmonic Blaschke sum} $K
\check{+}_{\varphi,\varepsilon}L$ of star bodies $K$ and $L$ in
${\Bbb R}^{n}$, which is an extension of the harmonic Blaschke
addition $\check{+}$, defined by (see section 3)
$$\varphi\left(\frac{k_{1}\rho(K,x)^{n}}{k\rho(K\check{+}_{\varphi,\varepsilon}L,x)^{n}}
,\frac{k_{2}\rho(L,x)^{n}}{k\rho(K\check{+}_{\varphi,\varepsilon}L,x)^{n}}\right)=1,
\eqno(1.5)$$ where $\varphi\in \Phi_{2}$ and $k_{1}, k_{2}, k$ are
positive constants.

If we taking for
$\varphi(x_{1},x_{2})=x^{(n+1)/n}_{1}+x^{(n+1)/n}_{2}$,
$k_{1}=1/V(K)^{n/(n+1)}, k_{2}=1/V(L)^{n/(n+1)}$ and
$k=1/V(K\check{+}_{\varphi,\varepsilon}L)^{n/(n+1)}$ in (1.5),
then the {\it Orlicz harmonic Blaschck addition} (1.5) reduces to
the harmonic Blaschck addition (1.3). On the other hand, for
$p\geq 1$, putting
$\varphi(x_{1},x_{2})=x_{1}^{(n+p)/n}+x_{2}^{(n+p)/n}$,
$k_{1}=V(K)^{-n/(n+p)}, k_{2}=V(L)^{-n/(n+p)}$ and
$k=V(K\check{+}_{\varphi,\varepsilon}L)^{-n/(n+p)}$ in (1.5), the
{\it Orlicz harmonic Blaschck addition} (1.5) reduces to the
following $L_{p}$-{\it harmonic Blaschck addition}, which was
defined in [5] and [46].
$$\frac{\rho(K\check{+}_{p}L,\cdot)^{n+p}}{V(K\check{+}_{p}L)}=\frac{\rho(K,\cdot)^{n+p}}
{V(K)}+\frac{\rho(L,\cdot)^{n+p}}{V(L)}.$$ Moreover, for $p\geq
n$, taking for $\varphi(x_{1},x_{2})=x_{1}^{p/n}+x_{2}^{p/n}$ and
$k_{1}=k_{2}=k=1$ in (1.5), the {\it Orlicz harmonic Blaschck
addition} (1.5) reduces to the $L_{p}$-{\it radial Minkowski
addition} (see e.g., [47]).
$$\rho(K\tilde{+}_{p}L,x)^{p}=\rho(K,x)^{p}+\rho(L,x)^{p},$$ where $K,L\in {\cal S}_{o}^{n}$ and $x\in {\Bbb R}^{n}$.

In section 5, we establish a new Brunn-Minkowski inequality for
{\it Orlicz harmonic Blaschke addition}.
$$1\geq\varphi\left(\frac{k_{1}V(K)}
{kV(K\check{+}_{\varphi,\varepsilon}L)},\frac{k_{2}V(L)}{kV(K\check{+}_
{\varepsilon,\varphi}L)}\right),\eqno(1.6)$$ with equality if and
only if $K$ and $L$ are dilates, where $\varphi\in \Phi_{2}$,
$K,L\in {\cal S}_{o}^{n}$, and $k_{1}, k_{2}, k$ are positive
constants. Taking for $\varphi(x_{1},x_{2})=x^{p}_{1}+x^{p}_{2}$,
$k_{1}=V(K)^{-1/p}, k_{1}=V(K)^{-1/p},$
$k=V(K\check{+}_{\varphi,\varepsilon}L)^{-1/p}$ and $p=(n+1)/n$ in
(1.6), (1.6) becomes to (1.4). For $p\geq 1$, putting
$\varphi(x_{1},x_{2})=x_{1}^{(n+p)/n}+x_{2}^{(n+p)/n}$,
$k_{1}=V(K)^{-n/(n+p)}, k_{2}=V(L)^{-n/(n+p)}$ and
$k=V(K\check{+}_{\varphi,\varepsilon}L)^{-n/(n+p)}$ in (1.6),
(1.6) becomes to
$$\frac{V(K\check{+}_{p} L)^{(n+p)/n}}{V(K\check{+}_{p} L)}\geq \frac{V(K)^{(n+p)/n}}{V(K)}+\frac{V(L)^{(n+p)/n}}{V(L)},\eqno(1.7)$$
with equality if and only if $K$ and $L$ are dilates. This is just
the $L_{p}$-dual Brunn-Minkowski inequality for $L_{p}$-harmonic
Blaschck addition (see e.g. [5]). Taking for $p=1$ in (1.7), (1.7)
reduces to (1.4). For $p\geq n$, taking for
$\varphi(x_{1},x_{2})=x_{1}^{p/n}+x_{2}^{p/n}$ and
$k_{1}=k_{2}=k=1$ in (1.6), (1.6) becomes to
$$V(K\tilde{+}_{p}L)^{p/n}\geq V(K)^{p/n}+V(L)^{p/n},$$
with equality if and only if $K$ and $L$ are dilates. This is just
the $L_{p}$-dual Brunn-Minkowski inequality (see [47]). For
different variants of the classical Brunn-Minkowski inequalities
we refer to [1], [2], [4], [38] and [42] and the references
therein.

In section 4, we consider the dual case of the Orlicz mixed
volumes and prove that there exist a new volume, {\it Orlicz dual
mixed volumes}, $\tilde{V}_{\varphi}(K,L)$, of star bodies $K$ and
$L$ in ${\Bbb R}^{n}$, by
$$\tilde{V}_{\varphi}(K,L)=:\frac{((\varphi)'_{l})(1)}{k_{1}}\lim_{\varepsilon\rightarrow 0^+}
\frac{kV(K\check{+}_{\varphi,\varepsilon}
L)-k_{1}V(K)}{\varepsilon}=\frac{1}{n}\int_{S^{n-1}}\varphi
\left(\frac{k_{2}\rho(L,u)^{n}}{k_{1}\rho(K,u)^{n}}\right)\rho(K,u)^{n}dS(u),\eqno(1.8)$$
where $\varphi\in \Phi$ and $K,L\in {\cal S}_{o}^{n}$, and
$k,k_{1},k_{2}$ are positive constants.

In Section 5, we establish an Orlicz Minkowski inequality for the
{\it Orlicz dual mixed volumes}.
$$\tilde{V}_{\varphi}(K,L)\geq V(K)\varphi\left(\frac{k_{2}}{k_{1}}\frac{V(L)}{V(K)}\right)
,\eqno(1.9)$$ with equality if and only if $K$ and $L$ are
dilates, where $\varphi\in \Phi$ and $K,L\in {\cal S}_{o}^{n}$,
and $k_{1},k_{2}$ are positive constants. Let
$\varphi(t)=t^{(n+p)/n}$ and $k_{1}=k_{2}=1$ in (1.8), then the
Orlicz mixed volumes $\tilde{V}_{\varphi}(K,L)$ reduces to the
well-known $L_{p}$-dual mixed volumes $\tilde{V}_{-p}(L,K)$ (see
[25]).
$$\tilde{V}_{-p}(L,K)=\frac{1}{n}\int_{S^{n-1}}\rho(L,u)^{n+p}\rho(K,u)^{-p}dS(u).\eqno(1.10)$$
On the other hand, putting $\varphi(t)=t^{(n+p)/n}$ and
$k_{1}=k_{2}=1$ in (1.9), then (1.9) reduces to the well-known
$L_{p}$-dual Minkowski inequality established by Lutwak in [25].
If $K,L\in {\cal S}_{o}^{n}$ and $p\geq 1$, then
$$\tilde{V}_{-p}(L,K)^{n}\geq V(L)^{n+p} V(K)^{-p},\eqno(1.11)$$
with equality holds if and only if $K$ and $L$ are dilates.
$L_{p}$ dual mixed volumes and inequalities for them play a
central role in the rapidly evolving $L_{p}$ Brunn-Minkowski
theory (see, e.g. [9], [10], [11], [12], [22], [27] and [41]).
Moreover, putting $k_{1}=k_{2}$ and $K=L$ in (1.8), the Orlicz
mixed volumes $\tilde{V}_{\varphi}(K,K)$ changes to the usual
volume $V(K)$.

In 2012, B\"{o}r\"{o}czky, Lutwak, Yang and Zhang [4] conjecture
that for origin-symmetric convex bodies $K$ and $L$ in ${\Bbb
R}^n$ and $0\leq \lambda\leq 1$,
$$V((1-\lambda)\cdot K +_{o}\lambda\cdot L)\geq V(K)^{1-\lambda}V(L)^{\lambda}
,\eqno(1.12)$$ where $(1-\lambda)\cdot K +_{o}\lambda\cdot L$
denotes the log Minkowski combination (see section 6). In fact, in
[4], they proved (1.12) only when $n=2$ and $K,L$ are
origin-symmetric convex bodies, and note that while it is not true
for general convex bodies. Moreover, they also shown that (1.12),
for all $n$, is equivalent to the following {\it log-Minkowski
inequality}.
$$\int_{S^{n-1}}\log\left(\frac{h(L,u)}{h(K,u)}\right)h(K,u)dS(K,u)
\geq V(K)\log\left(\frac{V(L)}{V(K)}\right),\eqno(1.13)$$ where
$S(K,u)$ is the surface area measure of $K$.

In section 6, we establish a {\it log dul-Minkowski inequality}.
If $K,L\in{\cal S}_{o}^{n}$ and ${\rm int}K\supset L$, then
$$\frac{1}{n}\int_{S^{n-1}}\log\left(1-\frac{\rho(L,u)^{n}}{\rho(K,u)^{n}}\right)\rho(K,u)^{n}dS(u)
\leq V(K)\log\left(1-\frac{V(L)}{V(K)}\right),\eqno(1.14)$$ with
equality if and only if $K$ and $L$ are dilates.

In 2010, the Orlicz projection body ${\bf\Pi} _{\varphi}$ of $K$
defined by Lutwak, Yang and Zhang [34]
$$h({\bf\Pi} _{\varphi},u)=\inf\left\{\lambda>0: \frac{1}{nV(K)}\int_{S^{n-1}}
\varphi\left(\frac{|u\cdot\upsilon|}{\lambda
h(K,\upsilon)}\right)h(K,\upsilon)dS(K,\upsilon)\leq
1\right\},\eqno(1.15)$$ for $u\in S^{n-1}$, $K\in{\cal
K}^{n}_{oo}$, and where ${\cal K}^{n}_{oo}$ denotes those sets in
${\cal K}^{n}$ containing the origin in their interiors. In [8],
the definition (1.15) of the Orlicz projection body suggested
defining, by analogy,
$$\widehat{V}_{\varphi}(K,L)=\inf\left\{\lambda>0: \frac{1}{nV(K)}\int_{S^{n-1}}
\varphi\left(\frac{h(L,u)}{\lambda h(K,u)}\right)h(K,u)dS(K,u)\leq
1\right\}.\eqno(1.16)$$

In section 7, we define a {\it dual Orlicz projection body},
$\tilde{\widehat{V}}_{\varphi}(K,L)$, as follows.
$$\tilde{\widehat{V}}_{\varphi}(K,L)=\inf\left\{\lambda>0: \frac{1}{nV(K)}\int_{S^{n-1}}
\varphi\left(\frac{k_{2}\rho(L,u)^{n}}{
k_{1}\lambda^{n}\rho(K,u)^{n}}\right)\rho(K,u)^{n}dS(u)\leq
1\right\},\eqno(1.17)$$ where $K,L\in {\cal S}_{o}^{n}$ and
$k_{1},k_{2}>0$. We prove an inequality for the {\it Orlicz dual
projection body}. Let $K,L\in{\cal S}_{o}^{n}$. If $\varphi\in
\Phi$ and $k_{1}, k_{2}>0$, then
$$\tilde{\widehat{V}}_{\varphi}(K,L)\geq\left(\frac{k_{2}}{k_{1}}
\left(\frac{V(L)}{V(K)}\right)\right)^{1/n}.\eqno(1.18)$$ If
$\varphi$ is strictly convex and $V(L)>0$, with equality if and
only if $K$ and $L$ are dilates. The inequality in special case
yields the $L_{p}$-dual Minkowski inequality (see Theorem 7.3 in
section 7).

\vskip 10pt \noindent{\large \bf 2 ~Definitions}\vskip 10pt

{\it 2.1 ~Dual mixed volumes}\vskip 10pt The radial Minkowski
linear combination,
$\lambda_{1}K_{1}\tilde{+}\cdots\tilde{+}\lambda_{r} K_{r},$
defined by
$$\lambda_{1}K_{1}\tilde{+}\cdots\tilde{+}\lambda_{r}
K_{r}=\{\lambda_{1}x_{1}\tilde{+}\cdots\tilde{+}\lambda_{r} x_{r}:
x_{i}\in K_{i},~ i=1,\ldots,r\}$$ for $K_{1},\ldots,K_{r}\in {\cal
S}_{o}^{n}$ and $\lambda_{1},\ldots,\lambda_{r}\in {\Bbb R}$. It
has the following important property: $$\rho(\lambda K\tilde{+}\mu
L,\cdot)=\lambda \rho(K,\cdot)+\mu\rho(L,\cdot),\eqno(2.1)$$ for
$K, L\in {\cal S}_{o}^{n}$ and $\lambda, \mu\geq 0.$

For $K_{1},\ldots,K_{r}\in {\cal S}^n$ and
$\lambda_{1},\ldots,\lambda_{r}\geq 0$, the volume of the radial
Minkowski linear combination
$\lambda_{1}K_{1}\tilde{+}\cdots\tilde{+}\lambda_{r}K_{r}$ is a
homogeneous $n$-th polynomial in the $\lambda_{i}$,
$$V(\lambda_{1}K_{1}\tilde{+}\cdots\tilde{+}\lambda_{r}K_{r})=\sum
\tilde{V}_{i_{1},\ldots,i_{n}}\lambda_{i_{1}}\cdots\lambda_{i_{n}},\eqno(2.2)
$$ where the sum is taken over all $n$-tuples
$(i_{1},\ldots,i_{n})$ whose entries are positive integers not
exceeding $r$. If we require the coefficients of the polynomial in
(2.2) to be symmetric in their argument, then they are uniquely
determined. The coefficient $\tilde{V}_{i_{1},\ldots,i_{n}}$ is
nonnegative and depends only on the bodies
$K_{i_{1}},\ldots,K_{i_{n}}$. Here we denote
$\tilde{V}_{i_{1},\ldots,i_{n}}$ to
$\tilde{V}(K_{i_{1}},\ldots,K_{i_{n}})$ and is called the dual
mixed volume of $K_{i_{1}},\ldots,K_{i_{n}}.$ If
$K_{1}=\cdots=K_{n-i}=K,$ $K_{n-i+1}=\cdots=K_{n}=L$, the dual
mixed volume is written as $\tilde{V}_{i}(K,L)$.

If $K_{1},\ldots,K_{n}\in {\cal S}_{o}^{n}$, the dual mixed volume
$\tilde{V}(K_{1},\ldots,K_{n})$ defined by ([28])
$$\tilde{V}(K_{1},\ldots,K_{n})=\frac{1}{n}\int_{S^{n-1}}\rho(K_{1},u)\cdots\rho(K_{n},u)dS(u).$$
For $K,L\in {\cal S}_{o}^{n}$, the $i$-th dual mixed volume of $K$
and $L$, $\tilde{V}_{i}(K,L)$, defined by
$$\tilde{V}_{i}(K,L)=\frac{1}{n}\int_{S^{n-1}}\rho(K,u)^{n-i}\rho(L,u)^{i}dS(u).\eqno(2.3)$$

\vskip 8pt {\it 2.2 ~$L_{p}$-dual mixed volumes}\vskip 10pt

For $p\geq 1$ and $K,L\in {\cal S}_{o}^{n}$, the $L_{p}$-dual
mixed volume, $\tilde{V}_{-p}(K,L)$, was defined in [25] by
$$-\frac{n}{p}\tilde{V}_{-p}(K,L)=\lim_{\varepsilon\rightarrow 0^{+}}\frac{V(K\hat{+}_{p}
\varepsilon\diamond L)-V(K)}{\varepsilon},\eqno(2.4)$$ where,
$K\hat{+}_{p}\varepsilon\diamond L$ is the $L_{p}$-harmonic
combination, is defined by
$$\rho(K\hat{+}_{p}\varepsilon\diamond L,\cdot)^{-p}=\rho(K,\cdot)^{-p}+\varepsilon \rho(L,\cdot)^{-p}.$$
The $L_{p}$-dual mixed volume $\tilde{V}_{-p}(K,L)$, for $K,L\in
{\cal S}_{o}^{n}$, was obtained by (see [25])
$$\tilde{V}_{-p}(K,L)=-\frac{p}{n}\lim_{\varepsilon\rightarrow 0^{+}}\frac{V(K\hat{+}_{p}\varepsilon
\diamond L)-V(K)}{\varepsilon}$$
$$~~~~~~~~~~~~~~~~~~=\frac{1}{n}\int_{S^{n-1}}\rho(K,u)^{n+p}\rho(L,u)^{-p}dS(u).\eqno(2.5)$$

\vskip 8pt {\it 2.3 ~Orlicz dual mixed volumes}\vskip 10pt

Gardner, Hug and Weil [8] defined the {\it Orlicz mixed volumes}
of convex bodies $K$ and $L$ as follows.
$$V_{\varphi}(K,L)=\frac{1}{n}\int_{S^{n-1}}\varphi\left(\frac{h(L,u)}{h(K,u)}\right)h(K,u)dS(K,u),$$
for all $K,L\in {\cal K}^{n}$ and $\varphi\in \Phi$, where
$S(K,u)$ is the surface area measure of $K$.

The {\it Orlicz dual mixed volumes}, $\tilde{V}_{\varphi}(K,L)$,
of star bodies $K$ and $L$ in ${\Bbb R}^{n}$, defined in section
4, by
$$\tilde{V}_{\varphi}(K,L)=:\frac{1}{n}\int_{S^{n-1}}\varphi
\left(\frac{k_{2}\rho(L,u)^{n}}{k_{1}\rho(K,u)^{n}}\right)\rho(K,u)^{n}dS(u),\eqno(2.6)$$
where $\varphi\in \Phi$ and $K,L\in {\cal S}_{o}^{n}$, and
$k_{1},k_{2}$ are positive constants.

\vskip 10pt \noindent{\large \bf 3 ~Orlicz harmonic Blaschke
addition}\vskip 10pt

In this section, we define a new concepts: {\it Orlicz harmonic
Blaschke addition}. Let $m\geq 2$ and $K_{j}\in {\cal S}_{o}^{n}$,
$k>0, k_{j}>0$, $j=1,\ldots,m$, we define the {\it Orlicz harmonic
Blaschke addition} of $K_{1},\ldots,K_{m}$, denoted by
$\check{+}_{\varphi}(K_{1},\ldots,K_{m})$, is defined by
$$\rho(\check{+}_{\varphi}(K_{1},\ldots,K_{m}),x)^{n}=\inf\left\{\lambda>0:
\varphi\left(\frac{k_{1}\rho(K_{1},x)^{n}}{k\lambda},\ldots,\frac{k_{m}\rho(K_{m},x)^{n}
}{k\lambda}\right)\leq 1\right\},\eqno(3.1)$$ for $x\in{\Bbb
R}^{n}.$ Equivalently, the Orlicz harmonic Blaschke addition
$\check{+}_{\varphi}(K_{1},\ldots,K_{m})$ can be defined
implicitly (and uniquely) by
$$\varphi\left(\frac{k_{1}\rho(K_{1},x)^{n}}{k\rho(\check{+}_{\varphi}(K_{1},
\ldots,K_{m}),x)^{n}},\ldots,\frac{k_{m}\rho(K_{m},x)^{n}}{k\rho(\check{+}_{\varphi}(K_{1},
\ldots,K_{m}),x)^{n}}\right)=1,\eqno(3.2)$$ if
$k_{1}\rho(K_{1},x)^{n}+\cdots+k_{m}\rho(K_{m},x)^{n}>0$ and by
$\rho(\check{+}_{\varphi}(K_{1},\ldots,K_{m}),x)=0$ if
$\rho(K_{1},x)=\cdots=\rho(K_{m},x)=0$, for all $x\in {\Bbb
R}^{n}$.

An important special case is obtained when
$$\varphi(x_{1},\ldots,x_{m})=\sum_{j=1}^{m}\varphi_{j}(x_{j}),$$
for some fixed $\varphi_{j}\in \Phi$ such that
$\varphi_{1}(1)=\cdots=\varphi_{m}(1)=1$. We then write
$\check{+}_{\varphi}(K_{1},\ldots,K_{m})=K_{1}\check{+}_{\varphi}\cdots\check{+}_{\varphi}K_{m}.$
This means that
$K_{1}\check{+}_{\varphi}\cdots\check{+}_{\varphi}K_{m}$ is
defined either by
$$\rho(K_{1}\check{+}_{\varphi}\cdots\check{+}_{\varphi}K_{m},u)^{n}=
\inf\left\{\lambda>0:\sum_{j=1}^{m}\varphi_{j}\left(\frac{k_{j}\rho(K_{j},x)^{n}}{k\lambda}
\right)\leq 1\right\},$$ for all $x\in {\Bbb R}^{n}$ and
$k,k_{j}>0$ $ j= 1,\ldots,m$, or by the corresponding special case
of (3.2).

{\bf Theorem 3.1}~ {\it If $\varphi\in\Phi_{m}$, then Orlicz
harmonic Blaschke addition $\check{+}_{\varphi}: ({\cal
S}_{o}^{n})^{m}\rightarrow {\cal S}_{o}^{n}$ is continuous,
monotonic, $GL(n)$ covariant, has the identity property.}

{\it Proof}~ We first prove that $+_{\varphi}$ is continuous. To
see this, indeed, let $K_{ij}\in {\cal S}_{o}^{n},$ $i \in {\Bbb
N}\cup\{0\},$ $j = 1,\ldots,m,$ be such that $K_{ij}\rightarrow
K_{0j}$ as $i\rightarrow\infty$. Notice that $\varphi$ is
continuous, we have
$$\varphi\left(\frac{k_{1}\rho(K_{01},x)^{n}}{k\lim_{i\rightarrow\infty}\rho(\check{+}_{\varphi}
(K_{i1},\ldots,K_{im}),x)^{n}},\ldots,\frac{k_{m}\rho(K_{0m},x)^{n}}{k\lim_{i\rightarrow\infty}\rho(\check{+}_{\varphi}
(K_{i1}, \ldots,K_{im}),x)^{n}}\right)=1.$$ Hence
$$\lim_{i\rightarrow\infty}\rho(\check{+}_{\varphi}
(K_{i1},\ldots,K_{im}),x)=\rho(\check{+}_{\varphi}(K_{01},\ldots,K_{0m}),x).$$
Next, we prove that $+_{\varphi}$ is monotonic. To see this, let
$K_{j}\subset L_{j}$ , where $K_{j}, L_{j}\in {\cal S}_{o}^{n},
j=1,\ldots, m.$ Let $x\in {\Bbb R}^{n}.$ If
$\rho(K_{1},x)=\cdots=\rho(K_{m},x)=0$, then
$\rho(\check{+}_{\varphi}
(K_{1},\ldots,K_{m}),x)=0\leq\rho(\check{+}_{\varphi}
(L_{1},K_{2},\ldots,K_{m}),x)$. If
$k_{1}\rho(K_{1},x)+\cdots+k_{m}\rho(K_{m},x)>0$, then
$k_{1}\rho(L_{1},x)+k_{2}\rho(K_{2},x)+\cdots+k_{m}\rho(K_{m},x)>0$
and using (3.2), $K_{1}\subset L_{1}$ , and the fact that
$\varphi$ is increasing in the first variable, we obtain
$$\varphi\left(\frac{k_{1}\rho(L_{1},x)^{n}}{k\rho(\check{+}_{\varphi}(L_{1},K_{2},
\ldots,K_{m}),x)^{n}},\frac{k_{2}\rho(K_{2},x)^{n}}{k\rho(\check{+}_{\varphi}(L_{1},K_{2},
\ldots,K_{m}),x)^{n}},\ldots,\frac{k_{m}\rho(K_{m},x)^{n}}{k\rho(\check{+}
_{\varphi}(L_{1},K_{2},\ldots,K_{m}),x)^{n}}\right)~~~~~~$$
$$~~~~=1=\varphi\left(\frac{k_{1}\rho(K_{1},x)^{n}}{k\rho(\check{+}_{\varphi}(K_{1},
\ldots,K_{m}),x)^{n}},\frac{k_{2}\rho(K_{2},x)^{n}}{k\rho(\check{+}_{\varphi}(K_{1},
\ldots,K_{m}),x)^{n}},\ldots,\frac{k_{m}\rho(K_{m},x)^{n}}{k\rho(\check{+}
_{\varphi}(K_{1}, \ldots,K_{m}),x)^{n}}\right)$$
$$~\leq\varphi\left(\frac{k_{1}\rho(L_{1},x)^{n}}{k\rho(\check{+}_{\varphi}(K_{1},
\ldots,K_{m}),x)^{n}},\frac{k_{2}\rho(K_{2},x)^{n}}{k\rho(\check{+}_{\varphi}(K_{1},
\ldots,K_{m}),x)^{n}},\ldots,\frac{k_{m}\rho(K_{m},x)^{n}}{k\rho(\check{+}
_{\varphi}(K_{1}, \ldots,K_{m}),x)^{n}}\right),$$ which again
implies that $\rho(\check{+}_{\varphi}(K_{1},
\ldots,K_{m}),x)\leq\rho(\check{+}_{\varphi}(L_{1},K_{2}
\ldots,K_{m}),x).$ By repeating this argument for each of the
other $(m-1)$ variables, we obtain
$\rho(\check{+}_{\varphi}(K_{1},
\ldots,K_{m}),x)\leq\rho(\check{+}_{\varphi}(L_{1},\ldots,L_{m}),x)$.

In the following, we prove $GL(n)$ covariant. Notice that
$$\varphi\left(\frac{k_{1}\rho(AK_{1},x)^{n}}{k\rho(\check{+}_{\varphi}(AK_{1},
\ldots,AK_{m}),x)^{n}},\ldots,\frac{k_{m}\rho(AK_{m},x)^{n}}{k\rho(\check{+}
_{\varphi}(AK_{1}, \ldots,AK_{m}),x)^{n}}\right)=1.$$ Hence
$$\varphi\bigg(\frac{k_{1}\rho(K_{1},A^{-1}x)^{n}}{k\rho(A^{-1}(\check{+}_{\varphi}(AK_{1},
\ldots,AK_{m})),A^{-1}x)^{n}},\ldots,\frac{k_{m}\rho(K_{m},A^{-1}x)^{n}}{k(\rho(\check{+}
_{\varphi}(AK_{1}, \ldots,AK_{m})),A^{-1}x)^{n}}\bigg)=1.$$ Set
$A^{-1}x = y,$ then
$$\varphi\bigg(\frac{k_{1}\rho(K_{1},y)^{n}}{k\rho(A^{-1}(\check{+}_{\varphi}(AK_{1},
\ldots,AK_{m})),y)^{n}},\ldots,\frac{k_{m}\rho(K_{m},y)^{n}}{k(\rho(\check{+}
_{\varphi}(AK_{1}, \ldots,AK_{m})),y)^{n}}\bigg)=1.$$ Hence
$$A^{-1}(\check{+}_{\varphi}(AK_{1},
\ldots,AK_{m}))=\check{+}_{\varphi}(K_{1}, \ldots,K_{m}).$$ This
shows Orlicz harmonic Blaschke addition $\check{+}_{\varphi}$ is
$GL(n)$ covariant.

The identity property is obvious from (3.2). $\Box$

In the section, we also define the Orlicz harmonic Blaschke linear
combination on the case $m=2$.

{\bf Definition 3.2}~ (Orlicz harmonic Blaschke linear
combination) {\it Orlicz harmonic Blaschke linear combination
$\check{+}_{\varphi}(K,L,\alpha,\beta)$ for $K,L\in {\cal
S}_{o}^{n}$, $k_{1}, k_{2}, k>0$ and} $\alpha,\beta\geq 0$ ({\it
not both zero}), {\it defined by
$$\alpha\varphi_{1}\left(\frac{k_{1}\rho(K,x)^{n}}{k\rho(\check{+}_{\varphi}(K,L,\alpha,\beta),x)^{n}}\right)+
\beta\varphi_{2}\left(\frac{k_{2}\rho(L,x)^{n}}
{k\rho(\check{+}_{\varphi}(K,L,\alpha,\beta),x)^{n}}\right)=1,\eqno(3.3)$$
if $\alpha k_{1}\rho(K,x)^{n}+\beta k_{2}\rho(L,x)^{n}>0$, and by
$\rho(\check{+}_{\varphi}(K,L,\alpha,\beta),x)=0$ if $\alpha
k_{1}\rho(K,x)^{n}+\beta k_{2}\rho(L,x)^{n}=0$, for all} $x\in
{\Bbb R}^{n}$.

It is easy to verify that when
$\varphi_{1}(t)=\varphi_{2}(t)=t^{(n+p)/n}, k_{1}=V(K)^{-n/(n+p)},
k_{1}=V(L)^{-n/(n+p)}, k=V(K+_{\varphi,\varepsilon}L)^{-n/(n+p)}$,
then the {\it Orlicz harmonic Blaschke linear combination}
$\check{+}_{\varphi}(K,L,\alpha,\beta)$ changes to the
$L_{p}$-harmonic Blaschke linear combination $\alpha\cdot
K\check{+}_{p}\beta\cdot L.$ Moreover, if
$\varphi_{1}(t)=\varphi_{2}(t)=t^{(n+1)/n}, k_{1}=V(K)^{-n/(n+1)},
k_{1}=V(L)^{-n/(n+1)}, k=V(K+_{\varphi,\varepsilon}L)^{-n/(n+1)}$,
then the Orlicz harmonic Blaschke linear combination
$\check{+}_{\varphi}(K,L,\alpha,\beta)$ equals the harmonic
Blaschke linear combination $\alpha\cdot K\check{+}\beta\cdot L.$

Henceforth we shall write $K\check{+}_{\varphi,\varepsilon}L$
instead of $\check{+}_{\varphi}(K,L,1,\varepsilon)$, for
$\varepsilon\geq 0$, and assume throughout that this is defined by
(3.1), where $\alpha=1, \beta =\varepsilon$, and
$\varphi_{1},\varphi_{2}\in \Phi$. The particular instance of
interest corresponds to using (1.5) with $\varphi(x_{1},
x_{2})=\varphi_{1}(x_{1})+\varepsilon\varphi_{2}(x_{2})$ for
$\varepsilon>0$ and some $\varphi_{1},\varphi_{2}\in \Phi$.

\vskip 10pt \noindent{\large \bf 4 ~Orlicz dual mixed
volumes}\vskip 10pt

In this section, we define a new concepts: {\it Orlicz dual mixed
volumes}.

{\bf Lemma 4.1}~ {\it If $K,L\in {\cal S}_{o}^{n}$ and $k,k_{1}$
are positive constants, then
$$K\check{+}_{\varphi,\varepsilon}L\rightarrow\left(\frac{k_{1}}{k}\right)^{1/n}\cdot K$$ in the radial Hausdorff metric as} $\varepsilon\rightarrow
0^{+}$.

{\it Proof}~ From Theorem 3.1, it is easy that the {\it Orlicz
harmonic Blaschke addition} is monotonous. Hence, it follows from
(3.3) with $\alpha=1$ and $\beta=\varepsilon$ and the preceding
remarks that for $u\in S^{n-1}$
$$\rho(K,u)\leq\rho(K\check{+}_{\varphi,\varepsilon}L,u)\leq\rho(K\check{+}_{\varphi,1}L,u),$$
for all $\varepsilon\in (0,1].$ If $\rho(K,u)>0$, we conclude from
$$\varphi_{1}\left(\frac{k_{1}\rho(K,u)^{n}}{k\rho(K\check{+}_{\varphi,\varepsilon}L,u)^{n}}\right)
+\varepsilon\varphi_{2}\left(\frac{k_{2}\rho(L,u)^{n}}{k\rho(K\check{+}_{\varphi,\varepsilon}L,u)^{n}}\right)=1$$
that if $\varepsilon\rightarrow 0^{+}$ and if a subsequence of
$\{\rho(K\check{+}_{\varphi,\varepsilon}L,u): \varepsilon\in(0,
1]\}$ converges to a constant
$\mu\geq\left(\frac{k_{1}}{k}\right)^{1/n}\rho(K,u)$, then
$$\varphi_{1}\left(\frac{k_{1}\rho(K,u)^{n}}{k\mu^{n}}\right)=\varphi_{1}\left(\frac{k_{1}\rho(K,u)^{n}}{k\mu^{n}}\right)
+0\cdot\varphi_{2}\left(\frac{k_{2}\rho(L,u)^{n}}{k\mu^{n}}\right)=1=\varphi_{1}(1),$$
and $k\mu^{n}\geq k_{1}\rho(K,u)^{n}$. Therefore, in this case,
$$\rho(K\check{+}_{\varphi,\varepsilon}L,u)\rightarrow\left(\frac{k_{1}}{k}\right)^{1/n}\cdot
\rho(K,u)$$ as $\varepsilon\rightarrow 0^{+}$.

If $\rho(K,u)=0$ and $\rho(L,u)>0$ then (3.3) with $\alpha=1$ and
$\beta=\varepsilon$ implies that
$$\varphi_{2}\left(\frac{k_{2}\rho(L,u)^{n}}
{k\rho(K\check{+}_{\varphi,\varepsilon}L,u)^{n}}
\right)=\frac{1}{\varepsilon}.$$ Hence
$$\rho(K\check{+}_{\varphi,\varepsilon}L,u)^{n}
=\frac{k_{2}}{k}\left(\varphi_{2}^{-1}\left(\frac{1}
{\varepsilon}\right)\right)^{-1}\rho(L,u)^{n},$$ and so
$$\rho(K\check{+}_{\varphi,\varepsilon}L,u)^{n}\rightarrow 0=\frac{k_{1}}{k}\cdot
\rho(K,u)^{n}$$ as $\varepsilon\rightarrow 0^{+}$. The latter also
holds if $\rho(K,u)=\rho(L,u)=0.$ ~~~~~~~~~~~~~~$\Box$

{\bf Theorem 4.2}~ {\it Let $\varphi\in \Phi_{2}$ and
$\varphi_{1}, \varphi_{2}\in \Phi$. If $K,L\in {\cal S}_{o}^{n}$
and $k,k_{1},k_{2}$ are positive numbers, then
$$\lim_{\varepsilon\rightarrow 0^+}\frac{k\rho(K\check{+}_{\varphi,\varepsilon}
L,u)^{n}-k_{1}\rho(K,u)^{n}}{\varepsilon}=\frac{k_{1}}{(\varphi_{1})_{l}'(1)}\cdot\varphi_{2}
\left(\frac{k_{2}\rho(L,u)^{n}}{k_{1}\rho(K,u)^{n}}\right)\rho(K,u)^{n}\eqno(4.1)$$
uniformly for} $u\in S^{n-1}.$

{\it Proof}~ From the hypotheses, we have for $\varepsilon>0$
$$\frac{k_{1}\rho(K,u)^{n}}{k\rho(K\check{+}_{\varphi,\varepsilon}
L,u)^{n}}=\varphi_{1}^{-1}\left(1-\varepsilon\varphi_{2}
\left(\frac{k_{2}\rho(L,u)^{n}}{k\rho(K\check{+}_{\varphi,\varepsilon}
L,u)^{n}}\right)\right).$$ Hence
$$\lim_{\varepsilon\rightarrow 0^+}\frac{k\rho(K\check{+}_{\varphi,\varepsilon}
L,u)^{n}-k_{1}\rho(K,u)^{n}}{\varepsilon}~~~~~~~~~~~~~~~~~~~~~~~~~~~~~~~~~~~~~~~~~~~~~~~~~~~~~~~~~~~~~~$$
$$=\lim_{\varepsilon\rightarrow 0^+}\frac{1}{\varepsilon}k\rho(K\check{+}_{\varphi,\varepsilon}
L,u)^{n}\left(1-\frac{k_{1}\rho(K,u)^{n}}{k\rho(K\check{+}_{\varphi,\varepsilon}
L,u)^{n}}\right)~~~~~~~~~~~~~~~~~~~~~~~~~~~~~~~$$
$$~~~~~~~~~~~~~~~~~~=\lim_{\varepsilon\rightarrow
0^+}k\rho(K\check{+}_{\varphi,\varepsilon}
L,u)^{n}\left(\varphi_{2}\left(\frac{k_{2}\rho(
L,u)^{n}}{k\rho(K\check{+}_{\varphi,\varepsilon}
L,u)^{n}}\right)\frac{1-\varphi_{1}^{-1}\left(1-\varepsilon\varphi_{2}
\left(\frac{k_{2}\rho(L,u)^{n}}{k\rho(K\check{+}_{\varphi,\varepsilon}
L,u)^{n}}\right)\right)}{1-\left(1-\varepsilon
\varphi_{2}\left(\frac{k_{2}\rho(L,u)^{n}}{k\rho(K\check{+}_{\varphi,\varepsilon}
L,u)^{n}}\right)\right)}\right). \eqno(4.2)$$ Let
$$z=\varphi_{1}^{-1}\left(1-\varepsilon\varphi_{2}
\left(\frac{k_{2}\rho(L,u)^{n}}{k\rho(K\check{+}_{\varphi,\varepsilon}
L,u)^{n}}\right)\right),\eqno(4.3)$$ and note that $z\rightarrow
1^{-}$ as $\varepsilon\rightarrow o^{+}.$ From (4.2), (4.3) and in
view of Lemma 4.1, we obtain
$$\lim_{\varepsilon\rightarrow 0^+}\frac{k\rho(K\check{+}_{\varphi,\varepsilon}
L,u)^{n}-k_{1}\rho(K,u)^{n}}{\varepsilon}~~~~~~~~~~~~~~~~~~~~~~~~~~~~~~~~~~~~~~~~~~~~~~~~~~~$$
$$=\lim_{\varepsilon\rightarrow
0^+}k\rho(K\check{+}_{\varphi,\varepsilon}
L,u)^{n}\cdot\varphi_{2}\left(\frac{k_{2}\rho(
L,u)^{n}}{k\rho(K\check{+}_{\varphi,\varepsilon}
L,u)^{n}}\right)\frac{1-z}{\varphi_{1}(1)-\varphi_{1}(z)}$$
$$=\frac{k_{1}}{(\varphi_{1})_{l}'(1)}\cdot\varphi_{2}\left(\frac{k_{2}\rho(L,u)^{n}}{k_{1}\rho(K,u)^{n}}\right)\rho(K,u)^{n}.~~~~~~~~~~~~~~~~~~~~~~~~~~~~~\eqno(4.4)$$
Moreover, the convergence is uniform for $u\in {\cal S}^{n-1}$.
Indeed, by (4.2), (4.3) and (4.4), it suffices to recall that by
Lemma 4.1
$$\lim_{\varepsilon\rightarrow
0^{+}}\rho(K\check{+}_{\varphi,\varepsilon}
L,u)=\left(\frac{k_{1}}{k}\right)^{1/n}\cdot\rho(K,u),$$ uniformly
on $S^{n-1}$. Hence
$$\lim_{\varepsilon\rightarrow 0^+}\frac{k\rho(K\check{+}_{\varphi,\varepsilon}
L,u)^{n}-k_{1}\rho(K,u)^{n}}{\varepsilon}=
\frac{k_{1}}{(\varphi_{1})_{l}'(1)}\cdot\varphi_{2}
\left(\frac{k_{2}\rho(L,u)^{n}}{k_{1}\rho(K,u)^{n}}\right)\rho(K,u)^{n}$$
uniformly for $u\in S^{n-1}.$ ~~~~~$\Box$

{\bf Theorem 4.3}~ {\it Let $\varphi\in \Phi_{2}$ and
$\varphi_{1}, \varphi_{2}\in \Phi$. If $K,L\in {\cal S}_{o}^{n}$
and $ k_{1}, k_{2}, k$ are positive constants, then}
$$\frac{((\varphi)'_{l})(1)}{k_{1}}\lim_{\varepsilon\rightarrow 0^+}
\frac{kV(K\check{+}_{\varphi,\varepsilon}
L)-k_{1}V(K)}{\varepsilon}=\frac{1}{n}\int_{S^{n-1}}\varphi_{2}
\left(\frac{k_{2}\rho(L,u)^{n}}{k_{1}\rho(K,u)^{n}}\right)\rho(K,u)^{n}
dS(u).\eqno(4.5)$$

{\it Proof}~ This follows immediately from Theorem 4.2 and
(2.3).~~~~~~~~$\Box$

Denoting by $\tilde{V}_{\varphi}(K,L)$, for any $\varphi\in\Phi$,
the integral on the right-hand side of (4.5) with $\varphi_{2}$
replaced by $\varphi$, we see that either side of the equation
(4.5) is equal to $\tilde{V}_{\varphi_{2}}(K,L)$ and hence this
new Orlicz dual mixed volume $\tilde{V}_{\varphi}(K,L)$ has been
born and therefore the new Orlicz dual mixed volume plays the same
role as $\tilde{V}_{-p}(K,L)$ in the $L_{p}$-dual Brunn-Minkowski
theory.

{\bf Definition 4.4}~ (Orlicz dual mixed volumes) {\it For
$\varphi\in \Phi$ and $k_{1}, k_{2}>0$, Orlicz dual mixed volumes,
$\tilde{V}_{\varphi}(K,L)$, defined by
$$\tilde{V}_{\varphi}(K,L)=:\frac{1}{n}\int_{S^{n-1}}\varphi
\left(\frac{k_{2}\rho(L,u)^{n}}{k_{1}\rho(K,u)^{n}}\right)\rho(K,u)^{n}dS(u),
\eqno(4.6)$$ for all} $K,L\in {\cal S}_{o}^{n}$.

\vskip 10pt \noindent{\large \bf 5 ~Orlicz dual Minkowski and
Brunn-Minkowski inequalities}\vskip 10pt

In the section, we need define a Borel measure in $S^{n-1}$,
$\tilde{V}_{n}(K,\upsilon),$ called as {\it dual normalized cone
measure}

{\bf Definition 5.1}~ {\it If $K\in {\cal S}_{o}^{n},$ {\it dual
normalized cone measure}, $\tilde{V}_{n}(K,\upsilon),$ defined by}
$$d\tilde{V}_{n}(K,\upsilon)=\frac{\rho(K,\upsilon)^{n}}{nV(K)}dS(\upsilon).
\eqno(5.1)$$

As dual case, the {\it normalized cone measure}, for convex body
$K$, $\bar{V}_{n}(K,\upsilon)$, by
$$d\bar{V}_{n}(K,\upsilon)=\frac{h(K,\upsilon)}{nV(K)}dS(K,\upsilon),\eqno(5.2)$$
which was defined by Gardner, Hug and Weil [8].

{\bf Lemma 5.2} [18] (Jensen's inequality) {\it Suppose that $\mu$
is a probability measure on a space $X$ and $g: X\rightarrow
I\subset {\Bbb R}$ is a $\mu$-integrable function, where $I$ is a
possibly infinite interval. If $\varphi: I\rightarrow {\Bbb R}$ is
a convex function, then
$$\int_{X}\varphi(g(x))d\mu(x)\geq\varphi\left(\int_{X}g(x)d\mu(x)\right).
\eqno(5.3)$$ If $\varphi$ is strictly convex, equality holds if
and only if $g(x)$ is constant for $\mu$-almost all} $x\in X$.

{\bf Lemma 5.3}~ {\it Let $0<a\leq\infty$ be an extended real
number, and let $I=[0,a)$ be a possibly infinite interval. Suppose
that $\varphi: I\rightarrow [0,\infty)$ is convex with
$\varphi(0)=0$. If $K,L\in{\cal S}_{o}^{n}$ are such that
$\sqrt[n]{k_{2}}L\subset{\rm int}(a\sqrt[n]{k_{1}}K)$, then
$$\frac{1}{nV(K)}\int_{S^{n-1}}\varphi
\left(\frac{k_{2}\rho(L,u)^{n}}{k_{1}\rho(K,u)^{n}}\right)\rho(K,u)^{n}dS(u)\geq
\varphi\left(\frac{k_{2}}{k_{1}}\frac{V(L)}{V(K)}\right)
.\eqno(5.4)$$ If $\varphi$ is strictly convex the equality holds
if and only if $K$ and $L$ are dilates.}

{\it Proof}~ By (2.3) with $i=0$ and $K=L$, the {\it dual
normalized cone measure} $\tilde{V}_{n}(K,u)$ is a probability
measure on $S^{n-1}$. If $\sqrt[n]{k_{2}}L\subset {\rm
int}(a\sqrt[n]{k_{1}}K)$, then
$$0\leq\frac{k_{2}\rho(L,u)^{n}}{k_{1}\rho(K,u)^{n}}<a$$ for all $u
\in S^{n-1}$. Therefore we can use Jensen's inequality (5.3), to
obtain
$$\frac{1}{nV(K)}\int_{S^{n-1}}
\varphi\left(\frac{k_{2}\rho(L,u)^{n}}{k_{1}\rho(K,u)^{n}}\right)\rho(K,u)^{n}dS(u)$$
$$=\int_{S^{n-1}}
\varphi\left(\frac{k_{2}\rho(L,u)^{n}}{k_{1}\rho(K,u)^{n}}\right)d\tilde{V}(u)$$
$$~~~~~~~~\geq\varphi\left(\frac{k_{2}}{nk_{1}V(K)}\int_{S^{n-1}}\rho(L,u)^{n}dS(u)\right)$$
$$=\varphi\left(\frac{k_{2}}{k_{1}}\frac{V(L)}{V(K)}\right).~~~~~~~~~~~~~~~~~~$$
In the following, we discuss the equal condition of (5.4). Suppose
the equality holds in (5.4) and since equality must hold in
Jensen's inequality (5.3) as well, when $\varphi$ is strictly
convex we can conclude from the equality condition for Jensen's
inequality that
$$\frac{1}{nV(K)}\int_{S^{n-1}}\frac{k_{2}\rho(L,u)^{n}}{k_{1}\rho(K,u)^{n}}\rho(K,u)^{n}dS(u)=\frac{k_{2}\rho(L,\upsilon)^{n}}{k_{1}\rho(K,\upsilon)^{n}},$$
for $S(\cdot)$-almost all $\upsilon\in S^{n-1}$. The shows that
the equality holds in (5.4) if and only if $K$ and $L$ are
dilates. ~~~~~~~~ $\Box$

{\bf Theorem 5.4}~ (Orlicz dual Minkowski inequality) {\it Let
$\varphi\in \Phi$. If $K,L\in {\cal S}_{o}^{n}$ and $k_{1},
k_{2}>0$, then
$$\tilde{V}_{\varphi}(K,L)\geq V(K)\varphi\left(\frac{k_{2}}{k_{1}}\frac{V(L)}{V(K)}\right)
,\eqno(5.5)$$ with equality if and only if $K$ and $L$ are
dilates.}

{\it Proof}~ This follows immediately from (4.6) and Lemma 5.3
with $a=\infty.$ ~~~~~~$\Box$

{\bf Corollary 5.5} ([25]) {\it If $K,L\in {\cal S}_{o}^{n}$ and
$p\geq 1$, then
$$\tilde{V}_{-p}(L,K)^{n}\geq V(L)^{n+p} V(K)^{-p},$$
with equality holds if and only if $K$ and $L$ are dilates.}

{\it Proof} The result follows immediately from Theorem 5.4 with
$\varphi(t)=t^{(n+p)/n}$ and $k_{1}=k_{2}=1$.~~~~~~$\Box$

{\bf Theorem 5.6}~ (Orlicz dual Brunn-Minkowski inequality) {\it
Let $\varphi\in \Phi_{2}$ and $k_{1}, k_{2}, k>0$. If $K,L\in
{\cal S}_{o}^{n}$, then}
$$1\geq\varphi\left(\frac{k_{1}V(K)}
{kV(K\check{+}_{\varphi,\varepsilon}L)},\frac{k_{2}V(L)}{kV(K\check{+}_{\varphi,\varepsilon}L)}
\right),\eqno(5.6)$$ {\it with equality if and only if $K$ and $L$
are dilates.}

{\it Proof}~ From the hypotheses and in view of the inequality
(5.5), we have
$$V(K\check{+}_{\varphi,\varepsilon}L)=\frac{1}{n}\int_{S^{n-1}}
\varphi\left(\frac{k_{1}\rho(K,u)^{n}}
{k\rho(K\check{+}_{\varphi,\varepsilon}
L,u)^{n}},\frac{k_{2}\rho(L,u)^{n}}{k\rho(K\check{+}_{\varphi,\varepsilon}
L,u)^{n}}\right)\rho(K\check{+}_{\varphi,\varepsilon}L)^{n}dS(u)~~~~~~$$
$$~~~~~~~~~~~~~~~~~~~~~~~~~=\frac{1}{n}\int_{S^{n-1}}\left(\varphi_{1}\left(\frac{k_{1}\rho(K,u)^{n}}
{k\rho(K\check{+}_{\varphi,\varepsilon}
L,u)^{n}}\right)+\varphi_{2}\left(\frac{k_{2}\rho(L,u)^{n}}
{k\rho(K\check{+}_{\varphi,\varepsilon}
L,u)^{n}}\right)\right)\rho(K\check{+}_{\varphi,\varepsilon}L)^{n}dS(u)$$
$$=\tilde{V}_{\varphi_{1}}(K\check{+}_{\varphi,\varepsilon}L,K)+
\tilde{V}_{\varphi_{2}}(K\check{+}_{\varphi,\varepsilon}L,L)~~~~~~~~~~~~~~~~~~~~~~~~~~~~~~~$$
$$\geq V(K\check{+}_{\varphi,\varepsilon}L)\left(\varphi_{1}
\left(\frac{k_{1}V(K)}{kV(K\check{+}_{\varphi,\varepsilon}L)}\right)+
\varphi_{2}\left(\frac{k_{2}V(L)}{kV(K\check{+}_{\varphi,\varepsilon}L)}\right)\right)~$$
$$=V(K\check{+}_{\varphi,\varepsilon}L)\varphi\left(\frac{k_{1}V(K)}
{kV(K\check{+}_{\varphi,\varepsilon}L)},\frac{k_{2}V(L)}{kV(K\check{+}_
{\varphi,\varepsilon}L)}\right).~~~~~~~~~~~~~~~~~$$

From the equality condition of Theorem 5.4, it follows that the
equality in (5.6) holds if and only if $K$ and $L$ are dilates.
~~~~$\Box$

{\bf Corollary 5.7} ([5]) {\it If $K,L\in {\cal S}_{o}^{n}$ and
$p\geq 1$, then
$$\frac{V(K\check{+}_{p} L)^{(n+p)/n}}{V(K\check{+}_{p} L)}\geq \frac{V(K)^{(n+p)/n}}{V(K)}+\frac{V(L)^{(n+p)/n}}{V(L)},$$
with equality if and only if $K$ and $L$ are dilates. }

{\it Proof} The result follows immediately from Theorem 5.6 with
$\varphi(x_{1},x_{2})=x_{1}^{(n+p)/n}+x_{2}^{(n+p)/n}$,
$k_{1}=V(K)^{-n/(n+p)}, k_{2}=V(L)^{-n/(n+p)}$ and
$k=V(K\check{+}_{\varphi,\varepsilon}L)^{-n/(n+p)}$.~~~~~~$\Box$

{\bf Corollary 5.8} ([47]) {\it If $K,L\in {\cal S}_{o}^{n}$ and
$p\geq n$, then
$$V(K\tilde{+}_{p}L)^{p/n}\geq V(K)^{p/n}+V(L)^{p/n},$$
with equality if and only if $K$ and $L$ are dilates.}

{\it Proof} The result follows immediately from Theorem 5.6 with
$\varphi(x_{1},x_{2})=x_{1}^{p/n}+x_{2}^{p/n}$ and
$k_{1}=k_{2}=k=1$.~~~~~~$\Box$

\vskip 10pt \noindent{\large \bf 6 ~Dual log-Minkowski
inequality}\vskip 10pt

Assume that $K,L$ is the nonempty compact convex subsets of ${\Bbb
R}^{n}$ containing the origin in their interiors, then the log
Minkowski combination, $(1-\lambda)\cdot K +_{o}\lambda\cdot L,$
defined by
$$(1-\lambda)\cdot K +_{o}\lambda\cdot L=\bigcap_{u\in S^{n-1}}
\{x\in {\Bbb R}^{n}:x\cdot u\leq
h(K,u)^{1-\lambda}h(L,u)^{\lambda}\},$$ for all real
$\lambda\in[0,1].$ B\"{o}r\"{o}czky, Lutwak, Yang and Zhang [4]
conjecture that for origin-symmetric convex bodies $K$ and $L$ in
${\Bbb R}^n$ and $0\leq \lambda\leq 1$,
$$V((1-\lambda)\cdot K +_{o}\lambda\cdot L)\geq V(K)^{1-\lambda}V(L)^
{\lambda}~ ?\eqno(6.1)$$ They [4] proved (6.1) only when $n=2$ and
$K,L$ are origin-symmetric convex bodies, and note that while it
is not true for general convex bodies. Moreover, they also shown
that (6.1), for all $n$, is equivalent to the following {\it
log-Minkowski inequality}
$$\int_{S^{n-1}}\log\left(\frac{h(L,u)}{h(K,u)}\right)d\bar{V}_{n}(K,\upsilon)
\geq\frac{1}{n}\log\left(\frac{V(L)}{V(K)}\right),\eqno(6.2)$$
where $\bar{V_{n}}(K,\cdot)$ is the {\it normalized cone measure}
for $K$. In fact, replacing $K$ and $L$ by $K+L$ and $K$,
respectively, (6.2) becomes to the following
$$\int_{S^{n-1}}\log\left(\frac{h(K,u)}{h(K+L,u)}\right)d\bar{V}_{n}(K+L,u)
\geq\log\left(\left(\frac{V(K)}{V(K+L)}\right)\right)^{\frac{1}{n}}.\eqno(6.3)$$
In [8], Gardner, Hug and Weil gave a log Minkowski type
inequality, for the nonempty compact convex subsets $K$ and $L$,
not origin-symmetric convex bodies, as follows. If $K$ is the
nonempty compact convex subsets of ${\Bbb R}^{n}$ containing the
origin in their interiors. $L$ is the nonempty compact convex
subsets of ${\Bbb R}^{n}$ containing the origin. Then
$$\int_{S^{n-1}}\log\left(\frac{h(K,u)}{h(K+L,u)}\right)d\bar{V}_{n}(K+L,u)
\leq\log\left(\frac{V(K+L)^{1/n}-V(L)^{1/n}}{V(K+L)^{1/n}}\right),\eqno(6.4)$$
with equality if and only if $K$ and $L$ are dilates or $L=\{o\}$.
They also shown that combining (6.3) and (6.4), may get the
classical Brunn-Minkowski inequality.
$$V(K+L)^{1/n}\geq V(K)^{1/n}+V(L)^{1/n},$$
whenever $K\in{\cal K}_{oo}^{n}$ and $L\in{\cal K}_{o}^{n}$ (where
${\cal K}_{o}^{n}$ is the class of members of ${\cal K}^{n}$
containing the origin) and (6.2) holds with $K$ and $L$ replaced
by $K+L$ and $K$, respectively. In particular, if (6.2) holds (as
it does, for origin-symmetric convex bodies when $n=2$), then
(6.2) and (6.4) together split the classical Brunn-Minkowski
inequality.

In this section, we establish a dual type log-Minkowski
inequalities.

{\bf Theorem 6.1}~ {\it Let $K, L\in{\cal S}_{o}^{n}$. If $k_{1}$
and $ k_{2}$ are positive numbers such that
$(\sqrt[n]{k_{2}}L)\subset{\rm int} (\sqrt[n]{k_{1}}K)$, then
$$\log\left(\frac{k_{1}V(K)-k_{2}V(L)}{k_{1}V(K)}\right)
\geq\int_{S^{n-1}}\log\left(\frac{k_{1}\rho(K,u)^{n}-k_{2}\rho(L,u)^{n}}{k_{1}
\rho(K,u)^{n}}\right)d\tilde{V}_{n}(u),\eqno(6.5)$$ with equality
if and only if $K$ and $L$ are dilates.}

{\it Proof}~ Since $K,L\in{\cal S}_{o}^{n}$ are such that
$(\sqrt[n]{k_{2}}L)\subset{\rm int} (\sqrt[n]{k_{1}}K)$. Let
$\varphi(t)=-\log(1-t)$, and notice that $\varphi(0)=0$ and
$\varphi$ is strictly increasing and strictly convex on $[0,1)$
with $\varphi(t)\rightarrow\infty$ as $t\rightarrow 1^{-}$. Hence
inequality (6.5) is a direct consequence of Lemma 5.3 with this
choice of $\varphi$ and $a=1$. ~~~~~~~~$\Box$

{\bf Corollary 6.2}~ {\it If $K, L\in{\cal S}_{o}^{n}$ such that
$L\subset {\rm int} K$, then
$$\log\left(\frac{V(K)-V(L)}{V(K)}\right)
\geq\int_{S^{n-1}}\log\left(\frac{\rho(K,u)^{n}-\rho(L,u)^{n}}{
\rho(K,u)^{n}}\right)d\tilde{V}_{n}(u),$$ with equality if and
only if $K$ and $L$ are dilates.}

This is just a dual form of the log Minkowski inequality (6.2).

{\bf Theorem 6.3}~ {\it Let $K,L\in{\cal S}_{o}^{n}$. If $k_{1},
k_{2}>0$ such that $(\sqrt[n]{k_{2}}L)\subset{\rm int}
(\sqrt[n]{k_{1}}(K\check{+}L))$, then
$$\log\left(\frac{k_{1}V(K\check{+}L)-k_{2}V(L)}{k_{1}V(K\check{+}L)}\right)
\geq\int_{S^{n-1}}\log\left(\frac{k_{1}\rho(K\check{+}L,u)^{n}-k_{2}\rho(L,u)^{n}}{k_{1}
\rho(K\check{+}L,u)^{n}}\right)d\tilde{V}_{n}(u),\eqno(6.6)$$ with
equality if and only if $K$ and $L$ are dilates.}

{\it Proof}~ If $K,L\in{\cal S}_{o}^{n}$, then
$K\check{+}L\in{\cal S}_{o}^{n}$. In view of
$(\sqrt[n]{k_{2}}L)\subset{\rm int}
(\sqrt[n]{k_{1}}(K\check{+}L))$ and from Theorem 6.1 with $K$
replaced by $K+L$, (6.6) easy follows.~~~~~~~~~~$\Box$

{\bf Corollary 6.4}~ {\it If $K,L\in{\cal S}_{o}^{n}$, then
$$\log\left(\frac{V(K\check{+}L)-V(L)}{V(K\check{+}L)}\right)
\geq\int_{S^{n-1}}\log\left(\frac{\rho(K\check{+}L,u)^{n}-\rho(L,u)^{n}}{
\rho(K\check{+}L,u)^{n}}\right)d\tilde{V}_{n}(u),$$ with equality
if and only if $K$ and $L$ are dilates.}

This is just a dual form of the log Minkowski type inequality
(6.4).

\vskip 10pt \noindent{\large \bf 7 Orlicz dual projection
body}\vskip 10pt

In 2010, the Orlicz projection body ${\bf\Pi} _{\varphi}$ of $K$
defined by Lutwak, Yang and Zhang [34]
$$h({\bf\Pi} _{\varphi},u)=\inf\left\{\lambda>0:\int_{S^{n-1}}
\varphi\left(\frac{|u\cdot\upsilon|}{\lambda
h(K,\upsilon)}\right)d\bar{V}_{n}(K,\upsilon)\leq
1\right\},\eqno(7.1)$$ for $K\in{\cal K}^{n}_{oo}, u\in S^{n-1},$
where $\bar{V_{n}}(K,\cdot)$ is the {\it normalized cone measure}
for $K$.

In [8], the definition (7.1) of the Orlicz projection body
suggested defining, by analogy,
$$\widehat{V}_{\varphi}(K,L)=\inf\left\{\lambda>0:\int_{S^{n-1}}
\varphi\left(\frac{h(L,u)}{\lambda
h(K,u)}\right)d\bar{V}_{n}(K,u)\leq 1\right\}.\eqno(7.2)$$

In this section, we define a dual of the Orlicz projection body,
as follows.

{\bf Definition 7.1}~ {\it If $K,L\in {\cal S}_{o}^{n}$ and
$k_{1},k_{2}>0$, then {\it Orlicz dual projection body},
$\tilde{\widehat{V}}_{\varphi}(K,L)$, defined by}
$$\tilde{\widehat{V}}_{\varphi}(K,L)=\inf\left\{\lambda>0:\int_{S^{n-1}}
\varphi\left(\frac{k_{2}\rho(L,u)^{n}}{
k_{1}\lambda^{n}\rho(K,u)^{n}}\right)d\tilde{V}_{n}(u)\leq
1\right\}.\eqno(7.3)$$

{\bf Theorem 7.2}~ {\it If $\varphi\in \Phi$ and $K, L\in{\cal
S}_{o}^{n}$ and $k_{1}, k_{2}>0$, then
$$\tilde{\widehat{V}}_{\varphi}(K,L)\geq\left(\frac{k_{2}}{k_{1}}
\left(\frac{V(L)}{V(K)}\right)\right)^{1/n}.\eqno(7.4)$$ If
$\varphi$ is strictly convex and $V(L)>0$, then equality holds if
and only if $K$ and $L$ are dilates.}

{\it Proof}~ Replacing $K$ by $\lambda K$, $\lambda>0$ in (5.4)
with $a=\infty$, we have
$$\int_{S^{n-1}}\varphi\left(\frac{k_{2}\rho(L,u)^{n}}{\lambda^{n}k_{1}
\rho(K,u)^{n}}\right)d\tilde{V}_{n}(u)
\geq\varphi\left(\frac{k_{2}}{\lambda^{n}k_{1}}\left(\frac{V(L)}{V(K)}\right)
\right).\eqno(7.5)$$ Let $\lambda>0$ such that
$$\int_{S^{n-1}}\varphi\left(\frac{k_{2}\rho(L,u)^{n}}{\lambda^{n}k_{1}
\rho(K,u)^{n}}\right)d\tilde{V}_{n}(u)\leq 1.$$ Hence
$$\varphi\left(\frac{k_{2}}{\lambda^{n}k_{1}}\left(\frac{V(L)}{V(K)}\right)
\right)\leq 1.$$ In view of $\varphi$ is strictly increasing in
each variable, we obtain
$$\left(\frac{k_{2}}{k_{1}}\left(\frac{V(L)}{V(K)}\right)\right)^{1/n}
\leq\lambda.\eqno(7.6)$$ From (7.3) and (7.6), (7.4) easy follows.

In the following, we discuss the equality condition of (7.4).
Suppose that equality holds, $\varphi$ is strictly convex and
$V(L)>0$. From (7.3), the exist
$\mu=\tilde{\widehat{V}}_{\varphi}(K,L)>0$ satisfies
$$\int_{S^{n-1}}\varphi\left(\frac{k_{2}\rho(L,u)^{n}}{\mu^{n} k_{1}
\rho(K,u)^{n}}\right)d\tilde{V}_{n}(u)=1.$$ Hence from the
assumption
$$\mu=\left(\frac{k_{2}}{k_{1}}\frac{V(L)}{V(K)}\right)^{1/n}.$$ Namely
$$\varphi\left(\frac{1}{\mu^{n}}\frac{k_{2}}{k_{1}}
\frac{V(L)}{V(K)}\right)=1.$$ Therefore the equality in (7.4)
holds for $\lambda=\mu$. But then the equality condition in (5.4)
shows that $\mu K$ and $L$ are dilates.~~~~~~~~~~~~$\Box$

{\bf Theorem 7.3}~ {\it The inequality} (7.4) {\it implies the
$L_{p}$-dual Brunn-Minkowski inequality.}

{\it Proof}~ Let $\lambda>0$ such that
$$\int_{S^{n-1}} \varphi\left(\frac{k_{2}\rho(L,u)^{n}}{
k_{1}\lambda^{n}\rho(K,u)^{n}}\right)d\tilde{V}_{n}(u)\leq
1.\eqno(7.7)$$ Putting $\varphi(t)=t^{(n+p)/n}, p\geq 1$ in (7.7),
we have
$$\int_{S^{n-1}}
\varphi\left(\frac{k_{2}\rho(L,u)^{n}}{
k_{1}\lambda^{n}\rho(K,u)^{n}}\right)d\tilde{V}_{n}(u)=\frac{1}{\lambda^{n+p}V(K)}\left(\frac{k_{2}}{k_{1}}\right)^{(n+p)/n}\tilde{V}_{-p}(L,K)\leq
1.\eqno(7.8)$$ From (7.3) and (7.8), as this case, we obtain
$$\tilde{\widehat{V}}_{\varphi}(K,L)=\left(\frac{k_{2}}{k_{1}}\right)^{1/n}\left(\frac{V_{-p}(L,K)}{V(K)}\right)^{1/(n+p)}.\eqno(7.9)$$
Putting (7.9) to (7.4), (7.4) reduces to
$$\tilde{V}_{-p}(L,K)^{n}\geq V(L)^{n+p}V(K)^{-p},$$
with equality if and only if $K$ and $L$ are dilates.

This is just the well-known $L_{p}$-dual Minkowski
inequality.~~~~~~~$\Box$

Moreover, there is an relationship between the Orlicz-Minkowski
inequalities (5.5) and (7.4). Since the two inequalities can be
written the following form, respectively
$$\frac{\tilde{V}_{\varphi}(K,L)}{ V(K)}\geq\varphi\left(\frac{k_{2}}{k_{1}}\frac{V(L)}{V(K)}\right)
,\eqno(7.10)$$
$$\varphi\left(\tilde{\widehat{V}}_{\varphi}(K,L)^{n}\right)\geq\varphi\left(\frac{k_{2}}{k_{1}}
\frac{V(L)}{V(K)}\right),\eqno(7.11)$$ Comparing (7.10) and
(7.11), we find that these are two inequalities of the different
precision. But when $\varphi(t)=t^{(n+p)/n}$, (7.10) and (7.11)
all change to the $L_{p}$-dual Minkowski inequality.

\end{document}